\documentclass[12pt]{article}
\usepackage{amsmath}
\usepackage{amssymb}
\usepackage{amscd}
\usepackage{amsthm}

\theoremstyle{plain}
\newtheorem{Thm}{Theorem}[section]

\theoremstyle{definition}
\newtheorem{Defn}[Thm]{Definition}
\newtheorem{Expl}[Thm]{Example}

\newtheorem{Que}[Thm]{Question}

\numberwithin{equation}{section}

\title{Finite generation of a canonical ring}
\author{Yujiro Kawamata}

\begin{document}

\maketitle

\begin{abstract}
The purpose of this note is to review an algebraic proof of the 
finite generation theorem due to 
Birkar-Cascini-Hacon-McKernan \cite{BCHM} whose method is based on the 
Minimal Model Program (MMP).
An analytic proof by Siu \cite{Siu-fg} will be reviewed by Mihai Paun.
\end{abstract}

\section{Introduction}

The finite generation of canonical rings
was a problem considered by Zariski \cite{Zariski}, and 
the proof in the case of dimension $2$ due to Mumford \cite{Mumford} 
in the appendix of \cite{Zariski}
is one of the motivations towards the minimal model theory of higher 
dimensional algebraic varieties.

Let $X$ be a smooth projective variety defined over a field $k$, 
$D$ a divisor on $X$, 
and $\mathcal{O}_X(D)$ the associated invertible sheaf.
Many problems in algebraic geometry are translated into questions on 
the vector space of holomorphic sections $H^0(X, D) = H^0(X,\mathcal{O}_X(D))$.
The Riemann-Roch problem is to determine this vector space. 
For example, the Riemann-Roch theorem tells us that the alternating sum
$\sum_{p=0}^n (-1)^p \dim H^p(X,\mathcal{O}_X(D))$ 
is expressed in terms of topological invariants.

Instead of considering a single vector space, we look at 
the graded ring 
\[
R(X,D) = \bigoplus_{m=0}^{\infty} H^0(X, mD)
\]
called the {\em section ring} for the pair $(X,D)$, where we use the additive 
notation for divisors instead of multiplicative one for sheaves; we have
$\mathcal{O}_X(mD) = \mathcal{O}_X(D)^{\otimes m}$.
There are obvious multiplication homomorphisms
\[
H^0(X, m_1D) \otimes H^0(X, m_2D) \to H^0(X, (m_1+m_2)D)
\]
and the section ring becomes a graded algebra over the base field 
$k = H^0(X, \mathcal{O}_X)$.
The following question arises naturally:

\begin{Que}
Is the section ring $R(X,D)$ finitely generated as a graded algebra 
over $k$?
\end{Que}

If $D$ is ample, then the answer is yes.
On the other hand, if $\dim X = 2$, 
then there exists an example where $R(X,D)$ is 
{\em not} finitely generated (\cite{Zariski}).
But the {\em canonical divisor} $K_X$, a divisor corresponding to 
the sheaf of holomorphic $n$-forms $\Omega^n_X$ for $n=\dim X$, has
a special status in algebraic geometry, and the
{\em canonical ring} $R(X, K_X)$ is finitely generated:

\begin{Thm}[Finite Generation Theorem]
Let $X$ be a smooth projective variety defined over a
field $k$ of characteristic $0$.
Then the canonical ring $R(X, K_X)$ is always finitely generated as a 
graded algebra over $k$.
\end{Thm}

The canonical divisor, or more precisely its linear equivalence class, is
special in many senses:

\begin{itemize}

\item It is naturally attached to any smooth variety as the determinant bundle 
of the cotangent bundle.

\item It is related to the Serre-Grothendieck duality.

\item The canonical ring is a {\em birational invariant}, 
i.e., if $f: X \to Y$ is a 
birational morphism between smooth projective varieties, then the natural
homomorphism $f_*: R(X,K_X) \to R(Y,K_Y)$ is an isomorphism.

\end{itemize}

There are two proofs, an algebraic one using the MMP 
by Birkar-Cascini-Hacon-McKernan 
\cite{BCHM} and an analytic one using complex analysis by Siu \cite{Siu-fg}.
We note that the algebraic proof yields 
the finite generation for varieties which are not necessarily of general type
because the proof is naturally \lq\lq logarithmic'' as we shall explain later. 

A geometric implication of the finite generation theorem is the 
existence of the {\em canonical model} $\text{Proj }R(X, K_X)$.
It is a traditional basic tool in the investigation 
of algebraic surfaces of general type.

The {\em anti-canonical ring} $R(X,-K_X)$ 
is not necessarily finitely generated.
For example, if $X$ is a ruled surface over a curve of genus greater than $1$,
then its anti-canonical ring is not finitely generated in general 
(\cite{Sakai}).
It is not a birational invariant either.

It is well known that canonical rings of algebraic curves 
are generated by elements of degree at most $3$, and those of 
algebraic surfaces are expected to be generated by elements of 
degree at most $5$.
But examples for varieties of dimension $3$ show that there is no bound
of degrees of generators for higher dimensional varieties even in
fixed dimensions.
A correct generalization for the effective statement is a birationality
question.
Hacon-McKernan \cite{HM-bounded} and Takayama \cite{Takayama-bounded} 
proved that 
there exists a number $m(n)$ depending only on the dimension $n$
such that the pluricanonical system $\vert m(n)K_X \vert$ gives 
a birational map for arbitrary $n$-dimensional variety $X$ of general type.
Chen-Chen \cite{CC} found an explicit and more realistic bound in the 
case of dimension three.
See also \cite{VZ} and \cite{Pacienza} for
development in this direction.

\vskip 1pc

The results of the four authors \cite{BCHM} include the complementary 
case of negative Kodaira dimension where the canonical rings tell very little.
The structure theorem as below, the existence of Mori fiber spaces, 
was another motivation for the MMP:

\begin{Thm}
Let $X$ be a smooth projective variety defined over a
field of characteristic $0$.
Assume that the canonical divisor $K_X$ is not pseudo-effective, i.e.,
not numerically equivalent to 
a limit of effective $\mathbf{Q}$-divisors.
Then there is a birational model of $X$ which has a Mori fiber 
space structure. 
\end{Thm}

A minimal model exists uniquely given a fixed birational 
class of surfaces of non-negative Kodaira dimension.
But the uniqueness of minimal models fails in higher dimensions even if
the Kodaira dimension is non-negative.
It is proved that arbitrary birationally equivalent minimal models 
are connected by
a sequence of operations called {\em flops}.
Birationally equivalent Mori fiber spaces are 
connected by {\em Sarkisov links} which are generalizations of 
elementary transformations of ruled surfaces.

\vskip 1pc

The third and more distant motivation for the MMP 
is a complete classification of algebraic varieties up to birational 
equivalence.
For example, if the Kodaira dimension is zero, 
then the canonical ring tells very little 
about the structure of the variety.
A minimal model of such a variety is still important and serves as a
starting point of more detailed structure theory.
We note that the existence of minimal models in general is still an open 
problem.

\vskip 1pc

We always assume that the charactersitic of the base field is zero 
in this paper.
We expect that the theorem is still true in positive characteristic, but
the known proofs depend heavily on characteristic zero methods. 
The restriction on the characteristic should be eventually removed and
the minimal model program (MMP) should work in arbitrary characteristic 
situation.
But this restriction is necessary for the proof presented here 
by the two reasons: 

\begin{itemize}

\item The resolution theorem of singularities by Hironaka is used extensively.
It would be possible to extend the resolution theorem over arbitrary 
characteristic base or even mixed characteristic base.

\item The Kodaira vanishing theorem holds true only in characteristic $0$.
There are counterexamples in positive characteristic.
The vanishing theorem is one of our main technical tools 
in every corner of the minimal model theory.  
In order to extend our theory to arbitrary characteristic, 
totally different technique should be developed.

\end{itemize}

We shall explain the following characteristic features of the MMP
(cf. \cite{KMM}):

\begin{itemize}

\item {\em inductive}. 
The basic idea is to use the induction on the dimension.
Then in more finer terms, the final product of the MMP, 
a minimal model or a Mori fiber space, 
is obtained by a step by step constructions 
called divisorial contractions and flips.

\item {\em logarithmic}. 
We consider pairs consisting of varieties and boundary 
divisors instead of varieties alone.

\item {\em relativistic}. 
The proof requires formulation in 
relative situation, i.e., we consider projective morphisms over base spaces.

\end{itemize}

We shall explain these points in the course of the algebraic proof.
It is important to note that 
the MMP is not yet completed because the following conjectures are still 
open:

\begin{itemize}

\item {\em The termination conjecture of flips}.
This is the conjecture which is still missing in order to prove the 
existence of a minimal model for arbitrary algebraic variety whose
canonical divisor is pseudo-effective.

\item {\em The Abundance Conjecture}.
We expect that some pluricanonical systems are base point free for 
arbitrary minimal models which are not necessarily of general type.
More precisely, if $f: (X,B) \to T$ is a projective morphism from a 
DLT pair such that $K_X+B$ is nef over $T$, 
then $K_X+B$ should be semi-ample over $T$,
i.e., numerically equivalent to a pull-back of an ample $\mathbf{R}$-divisor 
over $T$.

\end{itemize}

In particular, the following conjecture is still open: 
\lq\lq an algebraic variety with negative Kodaira dimension is uniruled''.

\vskip 1pc

The organization of this paper is as follows.
In order to motivate the use of the minimal models, 
we review a classical proof of the finite generation theorem 
in the case of dimension two in \S 2. 
We explain our point of view towrard minimal models, 
the numerical geometry, in \S3 and 
explain how complicated problems in birational geometry are 
simplified to those in linear algebra.
We explain how singularities of pairs appear naturally in \S4,
and try to justify these cumbersome definitions.
The algorithm of the Minimal Model Program (MMP) is explained in \S5.
The main results are presented in \S6 that there exists
a minimal model in the case where the 
canonical divisor is pseudo-effective or a Mori fiber space
otherwise under the additional assumption that the boundary is big.
The main idea of proof is to apply the vanishing theorems 
for the extension problems of pluricanonical forms (\S7).
The proofs of the main steps, the existence and termination of flips,
are respectively surveyed in \S8 and \S9. 
The last section \S10 is concerned on the non-uniqueness problem of the 
birational models.


\section{Case of dimension one or two: motivation of the MMP}

We start with the low dimensional cases in order to explain the idea how to
use the minimal models.
We note that the proofs of the finite generation theorem in dimensions one and
two are valid in arbitrary characteristic.

Let us consider the case $\dim X = 1$ and let $g$ be the genus of $X$.
There are three cases:

\vskip .5pc

{\em Case 1}: $g \ge 2$.  
Then $K_X$ is ample, 
hence the canonical ring is finitely generated.

\vskip .5pc

{\em Case 2}: $g=1$. 
Then $K_X = 0$, hence $R(X, K_X) \cong k[x]$.

\vskip .5pc

{\em Case 3}: $g = 0$. 
Then $R(X, K_X) \cong k$.

\vskip .5pc

Therefore the canonical ring is always finitely generated in dimension one.

\vskip 1pc

The above division into cases is generalized to higher dimensions by a
birational invariant called the {\em Kodaira dimension} of the pair 
$(X,D)$ defined by 
\[
\kappa(X,D) = \text{trans.deg }R(X,D) - 1
\]
if $R(X,D) \ne k$, and $\kappa(X,D) = -\infty$ otherwise.
The latter convention comes from the estimate that
$\dim H^0(X,mD) \sim m^{\kappa}$ for large $m$.
$\kappa(X,D)$ can take values among $-\infty, 0, 1, \dots , \dim X$.
In particular, we define $\kappa(X) = \kappa(X,K_X)$.
If $\kappa(X,D) = \dim X$, then $D$ is called {\em big}.
$X$ is said to be of {\em general type} if $K_X$ is big.

\vskip 1pc

We recall a proof by Mumford \cite{Mumford} (appendix to \cite{Zariski})
of the finite generation theorem when $\dim X = 2$ and $X$ is 
of general type.

\vskip .5pc

{\em Step 1}: A curve $C$ on $X$ is said to be a {\em $(-1)$-curve} if 
$C \cong \mathbf{P}^1$ and $(C^2) = -1$, where $(C^2)$ denotes the 
self intersection number.
In other words, the normal bundle $N_{C/X}$ has degree $-1$.
Castelnuovo's contraction theorem tells us that there exists a birational 
morphism $X \to X_1$ to another smooth projective surface 
which contracts $C$ to a point.
By using Castelnuovo's contraction theorem repeatedly, 
we obtain a birational morphism $X \to X'$ to a {\em minimal model} 
on which there is no $(-1)$-curve.
We note that the minimal model $X'$ is still smooth.
It is known that there exists unique minimal model $X'$ though the order of 
contractions of $(-1)$-curves can be arbitrary.

\vskip .5pc

{\em Step 2}: A curve $C$ on $X$ is said to be a {\em $(-2)$-curve} if 
$C \cong \mathbf{P}^1$ and $(C^2) = -2$.
Artin's contraction theorem (\cite{Artin}) states that all the $(-2)$-curves 
can be contracted to points by a birational morphism $X' \to X''$ to the 
{\em canonical model}.
The canonical model $X''$ is a projective surface 
with isolated rational singularities called Du Val singularities or
{\em canonical singularities}. 
Though $X''$ is singular, we have still define a canonical divisor 
$K_{X''}$ as a {\em Cartier} divisor, 
hence the canonical ring $R(X'',K_{X''})$.

\vskip .5pc

{\em Step 3}: The canonical ring is invariant under the contractions:
\[
R(X, K_X) \cong R(X',K_{X'}) \cong R(X'',K_{X''}).
\]
$K_{X''}$ is ample, hence $R(X,K_X)$ is finitely generated. 

\vskip .5pc

Step 1 will be generalized to the MMP, and Step 2 to the 
{\em Base Point Free Theorem}.

\vskip 1pc

A proof when $\dim X = 2$ and $\kappa = 1$ uses Kodaira's theory of 
elliptic surfaces (\cite{Kodaira-surface}).
 
\vskip .5pc

{\em Step 1} is the same as above.
Let $X'$ be a minimal model of $X$.
It is unique again.

\vskip .5pc

{\em Step 2}: By the Enriques classification, 
there is an elliptic surface structure $f: X' \to Y$; 
$Y$ is a smooth projective curve and general fibers of $f$ are elliptic curves. There may be degenerate fibers called {\em singular fibers}.

\vskip .5pc

{\em Step 3}: We define a {\em $\mathbf{Q}$-divisor}, a divisor with 
rational coefficients, $B_Y = \sum_i b_iB_i$ on $Y$ 
which measures the degeneration of the morphism $f$.
The coefficients $b_i$ belong to $\frac 1{12}\mathbf{Z}$, corresponding to the 
fact that there exists an automorphic form of degree $12$.
The values of the $b_i$ are  
determined by Kodaira's classification of singular fibers 
(\cite{Kodaira-surface}).
If a singular fiber of $f$ has finite local monodromy, then the corresponding
coefficient belong to the interval $(0,1)$.
Otherwise, the natural coefficient can be arbitrarily large, but we can 
move the divisor $B_Y$ to another $\mathbf{Q}$-divisor $B'_Y$ 
by using the automorphic form so that all 
the coefficients are less than $1$.
Thus we obtain the following {\em canonical bundle formula} of Kodaira:
\[
K_{X'} = f^*(K_Y+B'_Y).
\]
Therefore the canonical ring of $X$ 
is isomorphic to the {\em log canonical ring} for a lower dimensional pair 
$(Y,B'_Y)$: 
\[
R(X, K_X) \cong R(X',K_{X'}) \cong R(Y, K_Y+B'_Y).
\]
We note that the concept of the section ring is naturally extended to the 
case where the coefficients are no more integers as follows:
\[
R(X,D) = \bigoplus_{m=0}^{\infty} H^0(X, \llcorner mD \lrcorner)
\]
where $\llcorner mD \lrcorner$ denotes the {\em round down} of $mD$
obtained by taking the round downs of the coefficients.
Since $K_Y+B'_Y$ is ample, we have the finite generatedness 
of the canonical ring again.

This is one of the motivations of the \lq\lq log'' theory which will be 
explained later.

\vskip 1pc

The case $\dim X = 2$ and $\kappa=0$ is easy; we have
$R(X,K_X) \cong k[x]$.

The classification of surfaces tells us that a minimal model of 
such a variety, which is unique again, is isomorphic to
either an abelian surface, a K3 surface, or their quotient.
We note that the higher dimensional generalizations of this class 
include interesting 
varieties such as Calabi-Yau manifolds and hyperK\"ahler manifolds. 

\vskip 1pc

In order to explain the second motivation for the MMP, 
we consider the case where $\kappa(X) = - \infty$.
In this case, the canonical ring is trivial and tells very little.
Instead there is an explicit structure theorem when $\dim X = 2$ 
and $\kappa(X)=-\infty$.

\vskip .5pc

{\em Step 1} is the same.
Let $X'$ be a minimal model of $X$.

\vskip .5pc

{\em Step 2}: There are two cases according to the Enriques classification.
The first case is a ruled surface $f: X' \to Y$; 
$Y$ is a smooth projective curve and all fibers are isomorphic to 
$\mathbf{P}^1$.
In the second case, we have $X' \cong \mathbf{P}^2$. 

\vskip .5pc

Step 2 will be generalized to a {\em Mori fiber space} in the MMP
explained later.
This is another motivation for the MMP.

A minimal model $X'$ of an algebraic surface $X$ with $\kappa(X) \ge 0$
is uniquely determined by
the birational class of $X$ and is independent of the 
contraction process of $(-1)$-curves.
On the other hand, there are many minimal models in the sense that there 
are no $(-1)$-curves if $\kappa(X) = - \infty$.
It is known that birationally equivalent minimal models 
are connected by a sequence
of {\em elementary transformations}.

\vskip 1pc

There is another approach toward the finite generation theorem using the
{\em Zariski decomposition} (\cite{Zariski}).
In the case of a surface $X$, there is a birational morphism $f: X \to X'$
to its minimal model.
We can write $K_X = f^*K_X + E$, where $E$ is an effective divisor, a divisor
with non-negative coefficients, whose support coincides with the exceptional
locus, the locus where $f$ is not an isomorphism.
The point is that this decomposition can be constructed on $X$ without using
the minimal model provided that $\kappa(X) \ge 0$.

In general, if $D$ is an effective divisor on a smooth projective surface $X$,
then there exist uniquely determined $\mathbf{Q}$-divisors $P$ and $E$, 
divisors with coefficients in $\mathbf{Q}$, which satisfy the following 
conditions:

\begin{enumerate}

\item $P$ is {\em nef}, i.e., the intersection numbers $(P \cdot C)$ are 
non-negative for all curves $C$ on $X$. 

\item $E$ is effective.

\item If $\{E_1, \dots, E_t\}$ is the set of irreducible components of $E$,
then $(P \cdot E_i) = 0$ for all $i$ and the matrix $[(E_i \cdot E_j)]$ is
negative definite.

\end{enumerate}

For example, $K_X = f^*K_X + E$ is the Zariski decomposition of the canonical
divisor which is regarded as an effective $\mathbf{Q}$-divisor when
$mK_X$ is effective for a positive integer $m$.
Therefore a minimal model is virtually obtained as the nef part of the
canonical divisor.

We obtained a minimal model of a log surface by using
the Zariski decomposition in the paper \cite{log surface} 
as a log generalization
of the minimal model theory in dimension two.
Let $X$ be a smooth projective variety of dimension two and 
$B$ a reduced normal crossing divisor.
Then the Zariski decomposition $K_X+B=P+E$ gives a minimal model
of the pair $(X,B)$;
the support of $E$ can be contracted to points by a birational morphism
$f: X \to X'$ to a normal surface and we have $P = f^*K_{X'}$.
We note that the minimal model $X'$ may have singularities, 
called {\em log terminal singularities}, 
and that the coefficients of $P$ are not necessarily integers 
even if those of $B$ are equal to one.
Moreover we can prove that 
a linear system $\vert mP \vert$ 
is base point free for a positive integer $m$.

There is a higher dimensional generalization of the concept of the 
Zariski decomposition.
The technique of the base point free theorem explained later implies that the 
finite generation theorem holds if there exists a Zariski decomposition of a
(log) canonical divisor (\cite{ZD}, \cite{Moriwaki}).
But there is a counterexample for the existence of the Zariski decomposition 
if we consider arbitrary effective divisors (\cite{Nakayama-ZD}).
On the other hand, the Zariski decomposition of a log canonical divisor 
should exist because the log canonical divisor of a minimal model will give 
the Zariski decomposition.
The non-uniqueness of minimal models and the uniqueness of the 
Zariski decomposition are compatible, because birationally equivalent 
minimal models have equivalent canonical divisors.

We note that the Zariski decomposition can be achieved 
only on a blown-up variety of the original model if it ever exists; 
this is a phenomenon in dimension three 
or higher due to the existence of flips.
The point is to resolve simultaneously 
the base loci of all the pluricanonical systems on a fixed variety.
Such simultaneous resolution can be considered on the Zariski space,
the inverse limit of blowing up sequences.
The {\em $b$-divisors} by
Shokurov \cite{Shokurov-PL} is a formal concept to consider such situations.
The polytope decomposition theorem explained later will solve this 
infinity problem and lead to the flip theorem.

There is also a complex analytic version of the Zariski decomposition
by Tsuji (\cite{Tsuji-ZD}, \cite{DPS}).
It's existence is rather easily proved though it is not as strong as the 
algebraic counterpart.


\section{Numerical geometry}

We would like to generalize the minimal model theory 
to higher dimensional varieties.
The first task is to generalize contraction theorems of 
Castelnuovo and Artin. 
It is a highly non-linear problem to identify the locus to be contracted 
and to construct a contraction morphism.
A linearization of this problem is achieved by using the intersection 
numbers of divisors with curves.
The method is called {\em numerical geometry}.
Please distinguish it from enumerative geometry.
The idea to use cones in real vector spaces in order for 
the investigation of the birational geometry
goes back to Hironaka's thesis \cite{Hironaka}.

We consider a {\em relative} situation; 
let $X$ be an algebraic variety which is projective over a base space $T$.
If the base space $T$ is a point (absolute case), 
then $X$ is just a projective variety.
We shall need to consider more general base space $T$ 
in order to use inductive arguments.
Although the descriptions of the relative situations are slightly longer, 
the same methods of proofs work as in the absolute case.

We assume that $X$ is normal.
Let $D = \sum_i d_iD_i$ be an {\em $\mathbf{R}$-divisor},
a formal linear combination of prime divisors $D_i$, 
reduced irreducible subvarieties of codimension one, with coefficients $d_i$
in $\mathbf{R}$.
It is said to be an {\em $\mathbf{R}$-Cartier divisor} if it can be 
expressed as a linear combination of Cartier divisors with coefficients 
in $\mathbf{R}$.
As a counterpart, let $C = \sum_i c_iC_i$ be an {\em $\mathbf{R}$-$1$-cycle}, 
a formal linear combination of curves $C_i$, 
reduced irreducible subvarieties of 
dimension one, with coefficients $c_i$ in $\mathbf{R}$.
It is said to be {\em relative} for $f$ 
if the $C_i$ are mapped to points of the base space.
Then we can define a bilinear pairing called the {\em intersection number} 
$(D \cdot C)$ with values in $\mathbf{R}$ if $D$ is $\mathbf{R}$-Cartier and 
$C$ is relative for $f$.
Two $\mathbf{R}$-Cartier divisors or two relative $\mathbf{R}$-$1$-cycles 
are said to be {\em numerically equivalent} 
if they give the same intersection numbers 
when paired to arbitrary counterparts.
We note that we cannot define the intersection number 
of an arbitrary $\mathbf{R}$-divisor which is not $\mathbf{R}$-Cartier with 
a curve.

We consider real vector spaces 
\[
\begin{split}
&N^1(X/T) = \{\mathbf{R} \text{-Cartier divisor}\}/\text{numerical equivalence}
\\
&N_1(X/T) = \{\text{relative } \mathbf{R} \text{-} 1 \text{-cycle}\}
/\text{numerical equivalence}
\end{split}
\]
which are dual to each other and known to be finite dimensional.
The dimension 
$\rho(X/T) = \dim N^1(X/T)$ 
is called the {\em Picard number} of $X$ over $T$.
It is an important numerical invariant next to the dimension $\dim X$.

Let $\overline{NE}(X/T)$ be the closed convex cone in $N_1(X/T)$ 
generated by the numerical classes of relative curves, and let
$\text{Nef}(X/T)$ be the {\em nef cone} defined as 
the dual closed convex cone in $N^1(X/T)$:
\[
\text{Nef}(X/T) = \{v \in N^1(X/T) \,\vert\, 
(v \cdot C) \ge 0 \,\, \forall C\}
\] 
where the $C$ are arbitrary curves relative for $f$.
An $\mathbf{R}$-divisor is called {\em relatively nef}, {\em nef} over $T$, 
or {\em $f$-nef}, 
if its numerical class belongs to $\text{Nef}(X/T)$.
The following Kleiman's criterion is fundamental (\cite{Kleiman}):

\begin{Thm}
A Cartier divisor $D$ is relatively ample over $T$ 
if and only if its numerical class belongs to the interior 
$\text{Amp}(X/T)$ of 
the nef cone $\text{Nef}(X/T)$ called the {\em ample cone}.
\end{Thm}

An $\mathbf{R}$-Cartier divisor $D$ is said to be {\em ample} over $T$,
or {\em $f$-ample}, 
if its numerical class belongs
to the ample cone $\text{Amp}(X/T)$, i.e., if it is a linear 
combination of ample Cartier divisors with positive linear coefficients.
We note that an ample $\mathbf{R}$-Cartier divisor does not in general become 
a Cartier divisor by multiplying a positive number and 
it is not directly related to a projective embedding.

An $\mathbf{R}$-Cartier divisor $D$ is said to be {\em big} over $T$, 
or {\em $f$-big}, 
if one can write $D = A + E$ for an
$f$-ample $\mathbf{R}$-Cartier divisor $A$ and an {\em effective} 
$\mathbf{R}$-divisor $E$, an $\mathbf{R}$-divisor with non-negative 
coefficients.
In other words, an $\mathbf{R}$-divisor is big if it is bigger 
than an ample divisor.
The set $\text{Big}(X/T)$ of numerical classes of all the 
$f$-big $\mathbf{R}$-divisors is an open cone in $N^1(X/T)$ called the 
{\em big cone}.
$D$ is said to be {\em pseudo-effective} over $T$ 
if its numerical class belongs
to the the closure $\text{Psef}(X/T)$ of the big cone 
callled a {\em pseudo-effective cone}.
In other words, $D$ is pseudo-effective if its numerical class is a 
limit of effective $\mathbf{Q}$-divisors.
We have the following commutative diagram of inclusions:
\[
\begin{CD}
\text{Amp}(X/S) @>>> \text{Nef}(X/S) \\
@VVV @VVV \\
\text{Big}(X/S) @>>> \text{Psef}(X/S)
\end{CD}
\]
where the cones on the left are open and on the right their closures.

Let $f: X \to Y$ be an arbitrary surjective morphism  
to another normal variety which is projective over $T$ 
such that all the geometric fibers of $f$ are connected.
Then the natural homomorphism $f^*: N^1(Y/T) \to N^1(X/T)$ is 
injective, and the image of the nef cone satisfies
\[
f^*\text{Nef}(Y/T) = \text{Nef}(X/T) \cap f^*N^1(Y/T).
\]
If $f$ is not an isomorphism, then $F = f^*\text{Nef}(Y/T)$ 
is a face of $\text{Nef}(X/T)$
in the sense that $F \subset \partial \text{Nef}(X/T)$, because 
the pull-back of an ample divisor on $Y$ is never ample on $X$.
There is a dual face $F^*$ of $\overline{NE}(X/T)$ defined by
\[
F^* = \{w \in \overline{NE}(X/T) \,\vert\, f_*w = 0\}.
\]
Indeed we have 
\[
F^* = \{w \in \overline{NE}(X/T) \,\vert\, (D \cdot w) = 0 \forall D\}
\]
where the $D$ are arbitrary ample divisors on $Y$.
By the Zariski main theorem, the morphism $f$ is uniquely determined by the 
face $F$ or $F^*$ as long as $Y$ is normal.

The converse is not true; there are faces which do not correspond to
any morphisms.
The point of the MMP is that there always exists a morphism 
if the face $F^*$ on the cone of curves 
lies in the negative side for the log canonical divisor, i.e., 
if the intersection numbers with the log canonical divisor are negative
for all non-zero vectors in $F^*$.

In the case where $X$ is smooth and $T$ is a point, 
Mori proved in his famous paper \cite{Mori} 
that, if the canonical divisor $K_X$ is not nef, 
then there always exists an {\em extremal ray}, a face $F^*$ with 
$\dim F^* = 1$ on which $K_X$ is negative.
Moreover he proved that there exists a corresponding 
contraction morphism from $X$ in the case $\dim X = 3$.
These results were generalized to the cone and contraction theorems 
explained later. 

The key point in \cite{Mori} is to prove that there exists a rational curve
which generates an extremal ray $F^*$.
In order to prove this, he employed the deformation theory 
over a base field of positive characteristic, a purely algebraic method.
It is remarkable that this is still the only method to prove 
the existence of rational curves; 
there is no purely characteristic zero proof nor complex analytic proof.
The result is generalized to arbitrary pairs with log terminal singularities
in \cite{length} using \cite{MM}; 
an extremal ray is always generated by a rational curve. 

In this paper we shall only deal with cones in $N^1(X/T)$ 
and not consider those in $N_1(X/T)$.
The reason will be clear in the polytope decomposition theorem.


\section{Log canonical divisors and log terminal singularities}

In order to extend the concept of the minimality to higher dimensional 
varieties, 
we have to change the point of view.
We look at canonical divisors instead of varieties themselves.
Moreover, we have to deal with logarithmic pairs instead of varieties.
Singularities appear naturally in this way.

We explain how to compare canonical divisors of birationally equivalent 
varieties as a way to compare different birational models.
Let $\alpha: X \dashrightarrow X'$ be a birational map between 
smooth projective varieties.
Then there exists a third smooth projective variety $X''$ 
with birational morphisms $f:X''\to X$ 
and $f': X''\to X'$ such that $\alpha = f' \circ f^{-1}$.
If we define canonical divisors $K_X$ and $K_{X'}$ by using the 
same rational differential $n$-form, then the difference 
$f^*K_X - (f')^*K_{X'}$ is a well defined divisor which is supported 
above the locus where $\alpha$ is not an isomorphism.
Moreover, it is independent of the choice of the rational differential 
$n$-form which was used to define $K_X$ and $K_{X'}$.
We define that an equality $K_X \ge K_{X'}$ holds if and only if 
their pull-backs satisfy 
$f^*K_X \ge (f')^*K_{X'}$.

\vskip 1pc

Let us consider some examples before stating a formal definition of 
log terminal singularities.
If $f: X \to X'$ is a contraction morphism of a $(-1)$-curve $C$, 
then $K_X - C = f^*K_{X'}$.
An important observation is that the canoical divisor becomes smaller 
when the variety becomes smaller.
Therefore, we come to the following (temporary) defeinition:
{\em A variety $X$ is said to be a} minimal model
{\em if $K_X$ is minimal among all birationally equivalent varieties}.

We note that, if $f : X \to X'$ is a contraction morphism of 
$(-2)$-curves, then $K_X = f^*K_{X'}$ and $X$ has 
{\em canonical singularities}.

Let us consider an example in higher dimension.
Let $X$ be a smooth variety of dimension $n$ and $E$ a prime divisor 
which is isomorphic to $\mathbf{P}^{n-1}$ and such that the normal
bundle $N_{E/X}$ is isomorphic to $\mathcal{O}_{\mathbf{P}^{n-1}}(-d)$
for a positive integer $d$.
Then there exists a birational morphism $f: X \to X'$ which contracts
$E$ to a point.
If $d \ge 2$, then $X'$ has an isolated quotient singularity,
a typical example of a {\em log terminal singularity}.
We can still define a canonical divisor $K_{X'}$ as a 
{\em $\mathbf{Q}$-Cartier divisor}, namely $dK_{X'}$ is a Cartier
divisor.
We have an equation
\[
K_X+(1 - \frac nd)E = f^*K_{X'}.
\]
Thus $K_X > K_{X'}$ if $n > d$.
In this case, the singularity of $X'$ is said to be {\em terminal}.
For example, if $n=3$ and $d=2$, then $K_X > K_{X'}$.
This terminal singularity was first discovered by Mori \cite{Mori}. 

$d$ is the smallest positive integer for which $dK_{X'}$ becomes a Cartier 
divisor.  
Such an integer is called a {\em Cartier index} of the singularity.
Therefore there exists an isolated quotient singularity in dimension three 
which is terminal and has arbitrarily large Cartier index.
This is the reason why there is no bound of degrees for the generators
of canonical rings in dimension three.
Indeed let $X$ be a projective variety having only terminal singularities 
such that there is a point at which the Cartier index is $d$ and that the 
canonical divisor $K_X$ is ample.
Then the canonical ring of $X$ is {\em not} generated by elements whose degrees
are less than $d$.

\vskip 1pc

We have to consider the log version of a canonical divisor, 
i.e., the {\em log canonical divisor} $K_X+B$ of a log pair $(X,B)$ 
in the MMP.
There are many reasons which lead us to consider the log version:

\begin{itemize}

\item The canonical divisor of a fiber space is
described by using a log canonical divisor as explained already 
in the case of elliptic surfaces. 
The canonical divisor of a finite covering is similarly described. 

\item The adjunction formula for log canonical divisors is the basic tool
for the inductive argument on dimensions.  

\item The natural range for the power of vanishing theorems is 
the category of log varieties.

\end{itemize}

An objact of the {\em log theory} is a pair $(X,B)$ 
consisting of a normal variety $X$ 
and an $\mathbf{R}$-divisor $B$ called a {\em boundary}.
In many cases, we assume that $B$ is effective, i.e., the coefficients
are non-negative.
If $B$ is not effective, we often call $B$ a {\em subboundary} 
in order to distinguish from the effective case.
The coefficients of $B$ are not necessarily integers, but rational 
or even real numbers.
The rational coefficients naturally appear when we consider positive 
multiples of divisors, while
the real coefficients when taking the limits. 
The canonical divisor $K_X$ is still defined as a Weil divisor 
as long as the variety $X$ is normal.
The {\em log canonical divisor} is the sum $K_X+B$.
For example, if $X$ is smooth and $B$ is a {\em reduced} divisor, 
i.e., all the coefficients are equal to $1$, with normal crossing 
support, 
then $K_X+B$ is the divisor corresponding to the determinant line 
bundle of logarithmic differential forms $\Omega_X^1(\log B)$.
This is the origin of the name \lq\lq log''.

We note that our log structures are different from the log structures
of Fontaine-Illusie-Kato \cite{Kato}.
Both concepts are derived as 
generalizations of normal crossing divisors on smooth varieties.
There is no log theory which includes both at the moment.

We have to introduce singularities in higher dimensions by two resasons.
The first is that the canoical divisor can become smaller when the 
variety acquires singularities in the case $\dim X \ge 3$ as observed in 
the previous example.
The other is that we have to consider pairs of varieties and divisors, i.e., 
log varieties.

Now we state a formal definition:

\begin{Defn}
Let $X$ be a normal variety.
It is said to be {\em $\mathbf{Q}$-factorial} 
if an arbitrary prime
divisor on $X$ is a {\em $\mathbf{Q}$-Cartier divisor}, 
a $\mathbf{Q}$-linear combination of Cartier divisors.

Let $(X,B)$ be a pair of a normal variety $X$
and an $\mathbf{R}$-divisor $B = \sum_i b_iB_i$, and 
let $\mu: Y \to X$ be a projective birational morphism from another variety.
The {\em exceptional locus} of $\mu$ is the smallest closed subset of $Y$ 
such that $\mu$ is an isomorphism when restricted on its complement.
$\mu$ is said to be a {\em log resolution} (in a strict sense) 
of the pair $(X,B)$, if $Y$ is smooth, 
the exceptional locus is a normal crossing divisor, and moreover 
the union of the exceptional locus and the strict transforms $\mu_*^{-1}B_i$ 
of the prime divisors $B_i$ is a normal crossing divisor.

The pair $(X,B)$ is said to be {\em divisorially log terminal (DLT)} 
or to have only {\em divisorially log terminal singularities} 
if there exists a log resolution such that the 
following conditions are satisfied:

\begin{enumerate}

\item The coefficients $b_i$ belong to the interval $(0,1]$.

\item The log canonical divisor $K_X+B$ is an $\mathbf{R}$-Cartier divisor.

\item We can write 
\[
\mu^*(K_X+B) = K_Y + \mu_*^{-1}B + \sum d_jD_j
\]
with $d_j < 1$ for all $j$,
where $\mu_*^{-1}B = \sum b_i\mu_*^{-1}B_i$ is the strict transform of 
$B$ and the $D_j$ are prime divisors which are contained in the 
exceptional locus of $\mu$.

\end{enumerate}

If $0 < b_i < 1$ for all $i$ instead of the condition 1, 
then the pair is called {\em KLT}.
In other words, we can write
\[
\mu^*(K_X+B) = K_Y +B_Y
\]
where the coefficients of $B_Y$ belong the interval $(0,1)$.
The pair is said to be {\em log canonial (LC)} if all the coefficients of 
$B_Y$ belong to the interval $(0,1]$ in the above formula.

If $d_j \le 0$ or $d_j < 0$ for all $j$ instead of the condition 3, then
the pair is called {\em canonical} or {\em terminal}, respectively.
We usually assume that $B=0$ for canonical or terminal singularities,
but those with $B \ne 0$ are also useful.
\end{Defn}

For example, if $X$ is smooth and the support of $B$ is a normal crossing
divisor, then the pair $(X,B)$ is DLT (resp. KLT) if and only if 
the coefficients are in $(0,1]$ (resp. $(0,1)$).

The condition 2 is necessary for definig the pull-back in the condition 3.
The point of the log terminality is that $K_X+B$ is strictly 
smaller than $K_Y + \mu_*^{-1}B + \sum D_j$.
In other words, the pair is locally minimal among all birational models
in the log sense.
DLT and KLT pairs are respectively called {\em weal log terminal} and 
{\em log terminal} in \cite{KMM}.

The log canonical singularites do not behave 
so well as log terminal singularities.
For example, the underlying variety $X$ may have non Cohen-Macaulay 
singularities.
We can compare this with the fact that the underlying variety of a 
DLT pair has always rational singularities in characteristic zero.
But we have often a situation where there exists another boundary $B'$ on 
an LC pair $(X,B)$ such that $(X,B')$ is KLT.
In this case there is no problem because we can easily extend arguments 
in the MMP by usnig the perturbation to a new boundary 
$(1-\epsilon) B + \epsilon B'$.

To be KLT is an open condition under varying coefficients of the boundary,
LC is closed, and DLT is intermediate.
Each of them has advantage and disadvantage.
We can use a compactness argument for the LC pairs.

The KLT condition is equivalent to the $L^2$ condition in complex analysis.
This stronger version of the log terminality is easier to handle and 
more natural in some cases.
For example, KLT condition is independent of a log resolution.
Moreover, some statements are only true for KLT pairs.
On the other hand, it is necessary to consider 
general DLT pairs for the induction argument on the dimension.
Indeed the adjunction formula holds only along a boundary component
with coefficient one.

The pairs having $\mathbf{Q}$-factorial DLT singularities
can be characterized as those which may appear in the process of the MMP 
if we start with pairs of smooth varieties and reduced normal crossing 
divisors.
This is the category of pairs we work in the MMP.


\section{Minimal Model Program (MMP)}

We explain the algorithm of the MMP in this section.
There are some standard references including 
\cite{KMM}, \cite{FA}, \cite{KM}, \cite{Matsuki} and \cite{Corti-book}.

The cone and the contraction theorems are fundamental for the MMP.
As we already explained, some of the faces of the nef cone 
correspond to morphisms while some are not.
Mori's discovery in \cite{Mori} is that, if the dual face of the 
cone of curves is contained in the half space on which 
the canonical divisor is negative,
then there should always exist a corresponding morphism. 
In particular we define an {\em extremal ray} to be such 
a one dimensional face of the cone of curves.

Now we state the cone and contraction theorems for the log pairs (\cite{KMM}).
We state the theorem in the
dual terminology on the cone of divisors:

\begin{Thm}[Cone Theorem]
Let $f: (X,B) \to T$ be a projective morphism from a DLT pair, $H$ an ample 
divisor for $f$, and $\epsilon$ a positive number.
Then the part of the nef cone $\text{Nef}(X/T)$ which is visible from 
the point $[K_X+B+\epsilon H] \in N^1(X.T)$ is generated by 
finitely many points whose coordinates are rational numbers.
\end{Thm}

If $K_X+B$ is nef, then the assertion is empty.
If $K_X+B$ is not nef, then $K_X+B+\epsilon H$ is not nef either if 
$\epsilon$ is sufficiently small.
The Come Theorem says that the boundary of the nef cone facing the point 
$[K_X+B+\epsilon H]$ is like a boundary of a polyhedral cone, 
and consists of finitely many rational {\em faces}, 
intersections of $\text{Nef}(X/T)$ with linear subspaces of $N^1(X/T)$.
When the positive number $\epsilon$ approaches to zero, the number of 
visible faces increases, so that there may be eventually infinitely many faces.

\begin{Thm}[Contraction Theorem]
Let $f: (X,B) \to T$ be a projective morphism from a DLT pair, 
and let $F \subset \text{Nef}(X/T)$ be a face which is 
visible from the point $[K_X+B+\epsilon H]$ as above.
Then there exists a surjective morphism $\phi: X \to Y$ to a normal variety 
which is projective over the base space $T$ 
such that the geometric fibers of $\phi$ are connected and that 
\[
F = \phi^*\text{Nef}(Y/T).
\]
\end{Thm}

In particular, if the codimension of the face $F$ is equal to one, then
we call $\phi$ a {\em primitive contraction}.

The Contraction Theorem is equivalent to the Base Point Free Theorem:

\begin{Thm}[Base Point Free Theorem]
Let $f: (X,B) \to T$ be a projective morphism from a DLT pair,
and let $D$ be a Cartier divisor on $X$.
Assume that $D$ is $f$-nef and that $D - \epsilon(K_X+B)$ is $f$-ample 
for a small positive number $\epsilon$.
Then there exists a positive integer $m_0$ such that $mD$ is relatively 
base point free for arbitrary $m \ge m_0$.
\end{Thm}

The conclusion is saying that the natural homomorphism
\[
f^*f_*\mathcal{O}_X(mD) \to \mathcal{O}_X(mD)
\]
is surjective.
The Base Point Free Theorem implies in particular 
the {\em Non-Vanishing Theorem} 
which states that $f_*\mathcal{O}_X(mD) \ne 0$.

\vskip 1pc

We explain steps of the MMP formulated by Reid (\cite{Reid}) 
modified to the log and relative situation.
The coefficients of the boundaries are real numbers, so that the limiting
arguments are made possible.

We consider the category of {\em $\mathbf{Q}$-factorial DLT} 
pairs $(X, B)$ consisting of varieties $X$ which are 
{\em projective} over a base space $T$ together with 
$\mathbf{R}$-divisors $B$. 
Let $f: X \to T$ be the structure morphism.
This category includes pairs consisting of smooth varieties and 
normal crossing divisors with coefficients belonging to 
in the interval $(0,1]$.

\vskip .5pc

{\em Case 1}: If $K_X + B$ is relatively nef over the base (i.e., $f$-nef), 
then $(X, B)$ is {\em minimal}. 
We stop here in this case.

\vskip .5pc

{\em Case 2}: If $K_X + B$ is not relatively nef, 
then there exists a primitive contraction morphism
$\phi: X \to Y$ by the Cone and Contraction Theorems.
It is important to note that 
we need to choose one of the codimension one faces of the nef cone.

There are three cases for $\phi$.

\begin{itemize}

\item If $\dim X > \dim Y$, then $\phi$ is called a {\em Mori fiber space}.
We stop here in this case.

\item If $\phi$ is a birational morphism and contracts a divisor, 
i.e, if $\phi$ is  
a {\em divisorial contraction}, 
then the exceptional locus of $\phi$ is proved to be a prime divisor
thanks to the $\mathbf{Q}$-factoriality condition.
The new pair $(X',B') = (Y, \phi_*B)$ is again $\mathbf{Q}$-factorial and 
DLT.  
We go back to Case 1 or 2, and continue the program.

\item If $\phi$ is {\em small}, i.e., if $\phi$ changes the variety 
only in codimension two
or higher, then there exists another small projective birational morphism 
$\phi': X' \to Y$ for which $K_{X'} + B'$ is positive where $B'$ is the
strict transform of $B$.
We note that the set of prime divisors stay unchanged under the 
birational map $(\phi')^{-1} \circ \phi$ which is called a {\em flip}. 
The {\em existence of a flip}  
is now a theorem, called the {\em Flip Theorem}, by Hacon and McKernan 
\cite{HM-flip}.  
The new pair $(X', B')$ is again $\mathbf{Q}$-factorial and DLT.  
We can continue the program.

\end{itemize}

We should prove that there does not exist an infinite chain of flips
in order to obtain the final product, a minimal model or a Mori fiber space. 
This {\em Termination of Flips} is still a conjecture in the general case.

\vskip 1pc

It is important to note that an inequality $K_X + B > K_{X'} + B'$ always holds
after a divisorial contraction or a flip.

The Picard number drops by one after a divisorial contraction:
$\rho(Y/T) = \rho(X/T) - 1$.
Therefore the number of possible divisorial contractions is bounded 
by the Picard number.

But the Picard number stays the same after a flip.
Indeed there is a natural isomorphism $(\phi'_*)^{-1} \circ \phi_*:
N^1(X/T) \to N^1(X'/T)$.
Under the identification of real vector spaces by this isomorphism, 
the nef cones 
$\text{Nef}(X/T)$ and $\text{Nef}(X'/T)$ are adjacent
and touch each other along the face $F$, 
the pull-back of the nef cone $\text{Nef}(Y/T)$.
Therefore we can investigate the changes of birational models by 
looking at a partial cone decomposition of the real vector space.

\begin{Expl}
The following is an example of a flip called {\em Francia's flip}
\cite{Francia}.
This example was first considered as a counterexample to the existence
of a minimal model in higher dimension.
Now it is incorporated into the MMP.

We have $\dim X = 3$ and $B = 0$ in this example.
Let $X'$ be a smooth projective variety containing a 
{\em smooth rational curve} $C'$, i.e., a curve isomorphic to $\mathbf{P}^1$,
whose normal bundle is isomorphic to $\mathcal{O}_{\mathbf{P}^1}(-1) \oplus 
\mathcal{O}_{\mathbf{P}^1}(-2)$.
By blowing-up and blowing-down, we can construct 
a normal projective variety $X$  
with an isolated quotient singularity $P$ and a smooth rational curve 
$C$ on $X$ passing through $P$ such that $X' \setminus C'$ is isomorphic to 
$X \setminus C$.

There is a small contraction morphism $\phi: X \to Y$ 
which contracts $C$ to a point $Q$.
$K_X$ is negative for $\phi$; $(K_X \cdot C) = - \frac 12 < 0$.
There is a small biratiotnal morphism $\phi': X' \to Y$ on the other 
side such that $\phi'(C') = Q$.
$K_{X'}$ is positive for $\phi'$; $(K_{X'} \cdot C') = 1$.
\end{Expl}

The MMP is a way to obtain a minimal model. 
We note that this is not the only way.
We give a formal definition of a minimal model and a canonical model:

\begin{Defn}
Let $(X_0,B_0)$ be an LC pair which is projective over a base $T$.
We assume that there exists another boundary $B'_0$ such that 
$(X_0,B'_0)$ is KLT.
A new $\mathbf{Q}$-factorial LC pair $(X,B)$ 
projective over a base $T$ together with 
a birational map $\alpha: X_0 \dashrightarrow X$ over $T$ is said to be 
a {\em minimal model} of $(X_0,B_0)$ if the following conditions are satisfied:

\begin{enumerate}

\item $\alpha_0$ is surjective in codimension one, i.e., an arbitrary 
codimension one scheme theoretic point of $X$ is 
in the image of a morphism which represents $\alpha$.

\item $B = \alpha_*B_0$ is the strict transform. 

\item If $\mu_0: Y \to X_0$ and $\mu: Y \to X$ are common resolutions, then 
the difference $\mu_0^*(K_{X_0}+B_0) - \mu^*(K_X+B)$ is effective.
Moreover arbitrary codimension one point of $Y$ which remains codimension one 
on $X_0$ but not on $X$ is contained in
the support of the difference $\mu_0^*(K_{X_0}+B_0) - \mu^*(K_X+B)$.

\end{enumerate}

A surjective morphism $f: X \to Z$ with connected geometric fibers 
from a minimal model to a normal variety projective over $T$ is said to be  
the {\em canonical model} if $K_X+B$ is numerically equivalent to the pull-back
$f^*H$ for a relatively ample $\mathbf{R}$-divisor $H$ on $Z$.
\end{Defn}

The triple $(X,B,\alpha)$ is more precisely 
called a {\em marked minimal model}.
The last condition for the minimal model 
means that a prime divisor is contracted by $\alpha$ only
if there is a reason to be contracted.
The elliptic fibration considered in \S 2 is an example of a canonical model.

The Base Point Free Theorem implies that if the boundary $B$ of a minimal model $(X,B)$ is big over $T$, then a canonical model always exists.
The existence of a canonical model in general is a conjecture called the 
{\em Abundance Conjecture}.

A minimal model of a given pair may not be unique if it ever exists.
But they are always equivalent; if $(X,B)$ and $(X',B')$ are two minimal 
models of a pair $(X_0,B_0)$ over $T$, and if $\mu: Y \to X$ and 
$\mu': Y \to X'$ are common resolutions, then we have
$\mu^*(K_X+B) = (\mu')^*(K_{X'}+B')$.
It follows that a canonical model is unique if it exists.

\vskip 1pc

We shall need a modified version of the MMP, called the {\em MMP with scaling}
or the {\em directed MMP}.
It is much easier version than the general MMP.

The process of the MMP is not unique because we should choose a face 
in each step.
The scaled version of the MMP has smaller ambiguity of the choice.
The termination conjecture for arbitrary sequence of flips is in fact not 
necessary in order to prove the existence of a minimal model or
a Mori fiber space.
What we have to prove is that {\em some} sequence of flips terminates.
It turned out that the termination conjecture is easier
for the sequence of flips which is directed by an additional divisor.
We can expect that a sequence of flips has more tendency to terminate 
when the boundary moves along a line segment in the space of divisors.

We take an additional effective $\mathbf{R}$-Cartier divisor $H$ 
besides the pair $(X,B)$ such that $(X,B+H)$ is still DLT.
We assume that $K_X+B+H$ is $f$-nef.
Let 
\[
t_1 = \min \{t \in \mathbf{R}_{\ge 0} \,\vert\, K_X+B+tH \,\, 
\text{is }f\text{-nef}\} \in [0,1].
\]
If $K_X+B$ is not $f$-nef, 
then we choose a face $F$ such that $[K_X+B+t_1H] \in F$.
The existence of such a face is clear in the case where $B$ is big.
In the general case, the existence is proved by Birkar \cite{Birkar} 
using the boundedness of extremal rays \cite{length}.

We perform a divisorial contraction or a flip to the pair $(X,B)$.
Let us denote by $(X,B)$ and $H$ the new pair and the strict transform of $H$
by abuse of language.
Then the threshold $t_1$ for the new pair decreases or stays the same.
In other words, the MMP is directed by the {\em scaling} $H$.
If it reaches to $0$, then $K_X+B$ becomes $f$-nef, so we are done.

The point is that the log pair $(X,B+t_1H)$ is already minimal,
since $K_X+B+t_1H$ is $f$-nef, 
though it is an intermediate model of the MMP for $(X,B)$.
It will be proved, in the course of inductive argument, 
that the set of underlying varieties of minimal models  
is finite when we move the coefficients of the boundary on the line segment
joining $B$ and $B+H$ as above.
This is the Polytope Decomposition Theorem.
We shall apply this theorem to the MMP with scaling 
for the termination argument.


\section{Existence of minimal models}

The paper \cite{BCHM} proved that minimal models exist for a 
$\mathbf{Q}$-factorial DLT pairs $(X,B)$ 
in the case where the boundary divisor $B$ is big 
with respect to the pair $(X,B)$ over the base in the following sense:

\begin{Defn}
Let $(X,B)$ be an LC pair.
A subvariety $Z$ of $X$ is said to be an {\em LC center} of $(X,B)$ if 
there exists a log resolution $\mu: Y \to (X,B)$ such that, 
when one writes $\mu^*(K_X+B)=K_Y+B_Y$, there exists an irreducible
component $E$ of $B_Y$ whose coefficient is equal to one and such that 
$\mu(E)=Z$.
\end{Defn}

If $(X,B)$ is a DLT pair, then it is easy to see that $Z$ is an LC
center if and only if it is an irreducible component of the intersection of
some of the irreducible components of $B$ whose coefficients are equal to one.

\begin{Defn}
Let $f: (X,B) \to T$ be a projective morphism from a DLT pair.
An $\mathbf{R}$-Cartier divisor $D$ on $X$ is said to be 
{\em big with respect to $(X,B)$} if 
one can write $D = A + E$ for an ample $\mathbf{R}$-Cartier divisor $A$ 
and an effective $\mathbf{R}$-Cartier divisor $E$ whose 
support does not contain any LC center of $(X,B)$.

$D$ is called {\em pseudo-effective with respect to $(X,B)$} if $D+A$ is 
big with respect to $(X,B)$ for any ample $\mathbf{R}$-Cartier divisor $A$.
\end{Defn}

There are many advantages to consider big boundaries.
If $B$ is big with respect to $(X,B)$, then the following hold:

\begin{itemize}

\item The log canonical divisor $K_X+B$ is semi-ample, i.e.,
it is numerically equivalent to a pull-back of an ample 
$\mathbf{R}$-Cartier divisor by a morphism.
In other words, the Abundance Conjecture holds in this case.
In particular, if $B$ is a $\mathbf{Q}$-divisor, then 
there exists a positive integer $m$ such that 
$m(K_X+B)$ is relatively base point free over the base.
This is a consequence of the Base Point Free Theorem.

\item The number of faces on the nef cone which was considered in the MMP
is finite even if $\epsilon$ goes to zero.
This is a consequence of the Cone Theorem.

\item The number of {\em marked minimal models} for a fixed pair 
is finite, 
where a marked minimal model is a minimal model together 
with a birational map from the original pair. 
This is one of the results proved in \cite{BCHM}.

\end{itemize}

We note that there are examples where items 2 and 3 are false 
when $B$ is not big.

Hacon and McKernan \cite{HM-flip} proved the following 
existence theorem of the flip
(when combined with results in \cite{BCHM}):

\begin{Thm}[Flip Theorem]
Let $(X,B)$ be a $\mathbf{Q}$-factorial DLT pair which is projective 
over a base, and 
$\phi: X \to Y$ a small contraction morphism in the MMP for $K_X+B$. 
Then there exists a flip $\phi': X' \to Y$.
\end{Thm}

In order to prove the Flip Theorem, 
it is sufficient to consider the case where the pair $(X,B)$ is KLT and 
$B$ is a $\mathbf{Q}$-divisor
by perturbing the coefficients of $B$ 
because the ampleness is an open condition.
In this case, the existence of a flip is equivalent 
to the finite generatedness of the
sheaf of graded $\mathcal{O}_Y$-algebras 
\[
\bigoplus_{m=0}^{\infty} \phi_*\mathcal{O}_X(\llcorner m(K_X+B) \lrcorner)
\]
over $\mathcal{O}_Y$, which is a a special case of the Finite Generation 
Theorem.
In other words, the inductive procedure of the MMP decomposes
a difficult global finite generation problem into easier local 
finite generation problems.

The termination conjecture of flips has not yet been fully proved.
But \cite{BCHM} proved a partial termination theorem 
for directed flips with scaling: 

\begin{Thm}[Scaled Termination Theorem]
Let $f:(X,B)\to T$ be a projective morphism from a $\mathbf{Q}$-factorial 
DLT pair to a quasi-projective variety.  
Assume that $B$ is big over $T$ with respect to $(X,B)$.
Then the MMP for $(X,B)$ with scaling always 
terminates after a finite number of steps.
\end{Thm}

By the MMP, we have the following consequence:

\begin{Thm}[Existence of Models]\label{existence of models}
Let $f: (X,B) \to T$ be a projective morphism from a $\mathbf{Q}$-factorial 
DLT pair to a quasi-projective variety.
Assume that $B$ is big over $T$ with respect to $(X,B)$.
Then there exists is a birational map
$\alpha: (X,B) \dashrightarrow (X',B')$ over $T$
to a $\mathbf{Q}$-factorial DLT pair with a projective morphism
$f': X' \to T$ which satisfies the following conditions:

\begin{enumerate}

\item $\alpha$ is surjective in codimension one, and $B' = \alpha_*B$. 

\item If $\mu: Y \to X$ and $\mu': Y \to X'$ are common 
resolutions, then $\mu^*(K_X+B) - (\mu')^*(K_{X'}+B')$ is
an effective $\mathbf{R}$-divisor whose support contains strict
transforms of all the exceptional prime divisors of $\alpha$.

\item If $[K_X+B] \in \text{Psef}(X/T)$, then $(X',B')$ is a 
{\em minimal model} of $(X,B)$ over $T$, i.e., $K_{X'}+B'$ is $f'$-nef. 
On the other hand, if $[K_X+B] \not\in \text{Psef}(X/T)$, 
then $(X',B')$ has a Mori fiber space structure $\phi: X' \to Y$ over $T$,
where $\dim Y < \dim X'$ and $- (K_{X'}+B')$ is $\phi$-ample.
\end{enumerate}
\end{Thm}

Now we state the Finite Generation Theorem.
We note that we do not need to assume the bigness of the boundary:

\begin{Thm}[Finite Generation Theorem]
Let $f: (X,B) \to T$ be a projective morphism from a KLT 
pair to a quasi-projective variety.
Assume in addition that the boundary $B$ is a $\mathbf{Q}$-divisor.
Then the relative canonical ring
\[
R(X/T,K_X+B) = \bigoplus_{m=0}^{\infty} f_*\mathcal{O}_X
(\llcorner m(K_X+B) \lrcorner)
\]
is finitely generated as a graded $\mathcal{O}_T$-algebra.
\end{Thm}

We note that the assumption of $B$ being a $\mathbf{Q}$-divisor
is indispensable.
Indeed the statement is clearly 
false for a non-rational $\mathbf{R}$-divisor.
On the other hand, the limiting argument in the proof of the existence 
of minimal models makes it 
necessary to consider $\mathbf{R}$-divisors in the proof.

The proof of the Finite Generation Theorem is as follows.
If $B$ is big with respect to the pair, then the theorem follows from 
the Base Point Free theorem together with the existence theorem of a 
minimal model.
For the general case, we use the following reduction theorem (\cite{FM}):

\begin{Thm}
Let $f: (X,B) \to T$ be a projective morphism from a KLT pair 
to a quasi-projective variety.
Assume that the boundary $B$ is a $\mathbf{Q}$-divisor, and that
$\kappa(X_{\eta},(K_X+B) \vert_{X_{\eta}}) \ge 0$ for the 
generic fiber $X_{\eta}$ of $f$.
Then there exists a projective birational morphism $\mu: Y \to X$ from a 
smooth variety, 
a surjective morphism $f: Y \to Z$ to a smooth variety projective over $T$,
a $\mathbf{Q}$-divisor $C = \sum c_jC_j$ on $Z$ such that $\sum C_j$ 
is a normal crossing divisor and $c_i \in (0,1)$, and a 
$\mathbf{Q}$-divisor $L$ on $Z$ which is nef over $T$ such that
the following are satisfied:

\begin{enumerate}

\item $K_Z+C+L$ is big over $T$.

\item $f \circ \mu^{-1}$ induces an isomorphism of 
graded $\mathcal{O}_T$-algebras:
\[
R(X/T, K_X+B) \cong R(Z/T, K_Z+C+L).
\]
\end{enumerate}
\end{Thm}

The outline of the proof of this theorem is as follows. 
We consider a rational map, called the {\em Iitaka fibration}, 
defined by some 
positive multiple of $K_X+B$ whose image has relative dimension equal to 
$\kappa(X_{\eta},(K_X+B) \vert_{X_{\eta}})$ over $T$, 
and modify it by changing birational models and using
covering tricks.
Then we apply the semi-positivity theorem of the Hodge bundles in
\cite{sp}.

The rest of the proof of the Finite Generation Theorem is as follows.
We write $K_Z+C+L = A + E$ for an ample $\mathbf{Q}$-divisor $A$ and 
an effective $\mathbf{Q}$-divisor $E$.
Then we apply the big case to a new pair
$(1+\epsilon)(K_Z+C+L) = K_Z+(C+\epsilon E)+(L+\epsilon A)$
for a small positive rational number $\epsilon$. 

\vskip 1pc

The following Polytope Decomposition Theorem for varying boundaries 
was proved in \cite{BCHM} after Shokurov \cite{Shokurov-3fold-model}. 
This theorem is interesting in its own right besides its importance for the 
termination argument.
We formulate the theorem in a more exact form and 
in the case where the boundary is not necessarily big.
The first part is a decomposition according to the canonical models:

\begin{Thm}[Polytope Decomposition Theorem 1]
Let $(X,\bar B)$ be a $\mathbf{Q}$-factorial 
KLT pair with a projective morphism $f: X \to T$ to a 
base space, $B_1,\dots,B_r$ effective $\mathbf{Q}$-Cartier divisors, 
and $\bar V$ a polytope contained in the set
$\{B = \sum_i b_iB_i \,\vert\, b_i \in \mathbf{R}\} \cong \mathbf{R}^r$
such that the pairs $(X,B)$ are LC for all $B \in \bar V$.
Consider a closed convex subset
\[
V = \{B \in \bar V \,\vert\, [K_X+B] \in \text{Psef}(X/T)\}.
\]
Assume that for each $B \in V$, there exist a minimal model 
$\alpha: (X,B) \dashrightarrow (Y,C)$ and a canonical model
$g: Y \to Z$ for $f: (X,B) \to T$.
Moreover assume that there exists a real number $\epsilon > 0$ 
for each $B$ such that the morphism 
$g: (Y,\alpha_*B') \to Z$ for $B' \in \bar V$ 
has minimal and canonical models whenever $B' \in \text{Psef}(Y/Z)$ and
$\Vert B' - B \Vert \le \epsilon$, where $\Vert \Vert$ denotes 
the maximum norm of the coefficients.
Then there exists a finite decomposition to disjoint subsets 
\[
V = \coprod_{j=1}^s V_j
\]
with rational maps $\beta_j: X \dashrightarrow Z_j$
which satisfies the following conditions:

\begin{enumerate}
\item $B \in V_j$ if and only if $\beta_j$ gives the canonical model 
for $f: (X,B) \to T$.

\item The closures $\bar V_j$, hence $V$, are polytopes for all $j$.
Moreover, if $\bar V$ is a rational polytope, then
so are the $\bar V_j$ and $V$.
\end{enumerate}
\end{Thm}

The second part is a finer decomposition according to the minimal models:

\begin{Thm}[Polytope Decomposition Theorem 2]
Under the same assumptions as the above theorem, 
each $V_j$ are further decomposed into a finite disjoint union
\[
V_j = \coprod_{k=1}^t W_{j,k}
\]
which satisfies the following conditions: 
let $\alpha: X \dashrightarrow Y$ be a birational map such that
\[
W = \{B \in V \,\vert\, \alpha \text{ is a minimal model for } (X,B)\}
\]
is non-empty. Then

\begin{enumerate}
\item There exists an index $j$ such that $W \subset \bar V_j$.

\item If $W \cap V_j$ is non-empty for some $j$, 
then $W$ coincides with one of the $W_{j,k}$.

\item The closure $\bar W_{j,k}$ is a polytope for any $j$ and $k$.
Moreover, if $\bar V$ is a rational polytope, then
so are the $\bar W_{j,k}$.
\end{enumerate}
\end{Thm}

The existence of $\bar B$ is assumed only to ensure that the MMP works well.
In the case where the $B \in V$ are big with respect to the pairs,
the existence of minimal and canonical models
is already proved, and the conditions in the above theorems are always 
satisfied.

We note that the conclusions of the 
above theorems do not contradict with examples in \cite{CY}, where there 
are infinitely many chambers,
even if the boundaries are not big.
The reason is that our finite decomposition theorem is a statement 
of local nature in a sense.

The termination of the flips for the scaled MMP is an easy consequence of the 
Polytope Decomposition Theorem.
Indeed the models before and after a flip correspond to different 
chambers.
There are only finitely many chambers on the line segment on which 
the scaled MMP is played, hence the termination.


\section{The vanishing theorems and the extension techniques}

Now we start to explain some ideas of proofs of the theorems.
The vanishing theorems of Kodaira type have surprizingly diverse applications 
in the birational geometry.
We need to assume that the characteristic of the base field is zero.

We start with the original Kodaira vanishing theorem (\cite{Kodaira-vanish}):

\begin{Thm}[Kodaira Vanishing Theorem]
Let $X$ be a smooth projective variety, and $D$ a divisor.  
If $D - K_X$ is ample, 
then $H^p(X, D) = 0$ for $p > 0$.
\end{Thm}

By applying a covering trick to the Kodaira vanishing theorem, 
we deduce its log and relative version
(\cite{vanishing}, \cite{Viehweg}, \cite{KMM}):

\begin{Thm}[Kawamata-Viehweg Vanishing Theorem]\label{KV}
Let $X$ be a normal variety, $f: X \to T$ a projective morphism 
to a base, 
$B$ an effective $\mathbf{R}$-divisor, and $D$ a Cartier divisor.  
If $(X,B)$ is KLT and $D -(K_X + B)$ is $f$-ample, 
then $R^pf_*\mathcal{O}_X(D) = 0$ for $p > 0$.
\end{Thm}

This generalization is used as in the following way.
Let $D$ be a divisor on a smooth projective variety $X$.
Suppose that $D - K_X$ is not ample, but close to ample.
We can sometimes find a small effective 
$\mathbf{R}$-divisor $B$ by a {\em perturbation} such that $(X,B)$ is KLT and 
$D - (K_X+B)$ is ample.
Then we have still the vanishing $H^p(X, D) = 0$ for $p > 0$.

The Cone and Contraction Theorems are proved by using the above 
Vanishing Theorem (\cite{KMM}).
The method of proofs are called {\em $x$-method}, 
which is a cohomological technique 
for proving base point freeness using the vanishing theorem.
The name came from the surprizingly similar applications of the  
vanishing theorem toward apparently different first two problems 
in Reid's list (\cite{Reid}).
The same method was later applied toward Fujita's conjecture
on the base point freeness of adjoint systems. 

The idea of the proof of the base point free theorem is as follows.
The point is to extend a holomorphic section from a codimension one subvariety 
to the whole space, and use the induction on the dimension.
Suppose that we want to prove the base point freeness of a
complete linear system
$\vert mD \vert$ on a smooth projective variety $X$ for a large integer $m$.
First we take a positive integer $m_1$ such that $\vert m_1D \vert$ 
is non-empty, 
where the existence of such an $m_1$ is guaranteed by the 
{\em Non-Vanishing Theorem} which is also proved by the $x$-method.
We want to make the base locus smaller by taking a larger integer $m_2$ 
so that the base locus disappears eventually.
We construct a projective birational morphism $\mu: Y \to X$ from a 
smooth variety with an effective divisor $E$ and a smooth prime divisor $Z$
such that:

\begin{enumerate}

\item The support of $E$ is contracted by $\mu$.

\item The image of $Z$ by $\mu$ is contained in the base locus 
of $\vert m_1D \vert$.

\end{enumerate}

The argument for finding such a situation looked tricky when it was found,
but we now know that the {\em log canonical threshold} is the concept hidden in
this argument.
We consider an exact sequence
\[
\begin{split}
&0 \to \mathcal{O}_Y(m_2\mu^*D + E - Z) \to \mathcal{O}_Y(m_2\mu^*D + E) \\
&\to \mathcal{O}_Z((m_2\mu^*D + E) \vert_Z) \to 0
\end{split}
\]
for a multiple $m_2$ of $m_1$.
If $m_2$ is suitably large, then we prove that 
\[
H^0(Z, (m_2\mu^*D + E) \vert_Z) \ne 0
\]
by using the Non-Vanishing Theorem,
and 
\[
H^1(Y, m_2\mu^*D + E - Z) = 0
\]
by using the Vanishing Theorem.
By condition 1, the natural homomorphism
\[
\mu_*: H^0(Y, m_2\mu^*D + E) \to H^0(X,m_2D)
\]
is bijective.
It follows that $\mu(Z)$ is not contained in the base locus 
of $\vert m_2D \vert$, hence it is strictly 
smaller that that of $\vert m_1D \vert$.

The proof of the Non-Vanishing Theorem is similar,
but we take an artificial non-complete 
linear system compared to the natural one 
$\vert mD \vert$.
The proof of the Cone Theorem is also similar.

\vskip 1pc

We worked \lq\lq upstairs'' in the above proof, 
i.e., on a resolution which lies 
above the original variety, and extend holomorphic sections from a divisor to
the whole variety by using the Vanishing Theorem for line bundles. 
If we use the vanishing theorem for multiplier ideal sheaves explained below, 
we can work \lq\lq downstairs'', i.e., on the original variety,
for the extension argument.
Then we can obtain more powerful extension theorems 
because it becomes possible 
to consider an infinite series of linear systems simultaneously.

We extend the vanishing theoren to the case where a pair $(X,B)$ 
consisting of a normal variety and an effective 
$\mathbf{R}$-divisor such that $K_X+B$ is an $\mathbf{R}$-Cartier divisor but 
the pair is {\em not} necessarily KLT.

\begin{Defn}
Let $(X,B)$ be a pair consisting of a normal variety and an effective 
$\mathbf{R}$-divisor such that $K_X+B$ is an $\mathbf{R}$-Cartier divisor.
Let $\mu: Y \to X$ be a log resolution and we write 
$\mu^*(K_X+B) = K_Y + \sum e_jE_j$, where $\sum E_j$ is a normal crossing 
divisor.
The {\em multiplier ideal sheaf} $I(X,B) \subset \mathcal{O}_X$ is defined by 
the following formula 
\[
I(X,B) = \mu_*\mathcal{O}_Y(- \sum \llcorner e_j \lrcorner E_j).
\]
\end{Defn}

It is independent of the choice of the log resolution.
We have $I(X,B) = \mathcal{O}_X$ if and only if $e_j < 1$ for all $j$, 
i.e., if $(X,B)$ is KLT.

Let $f: X \to T$ be a projective morphism and $D$ a Cartier divisor
on $X$.
If $D - (K_X+B)$ is $f$-ample, then Theorem~\ref{KV} implies that 
\[
R^p\mu_*\mathcal{O}_X(\mu^*D - \sum \llcorner e_j \lrcorner E_j) = 0
\]
and 
\[
R^p(f \circ \mu)_*\mathcal{O}_X(\mu^*D - \sum \llcorner e_j \lrcorner E_j) = 0
\]
for $p > 0$.
Hence we have a vanishing theorem for non-KLT pairs:

\begin{Thm}[Nadel Vanishing Theorem]
Let $f: X \to T$ be a projective morphism from a normal variety, 
$B$ an effective $\mathbf{R}$-divisor 
such that $K_X+B$ is an $\mathbf{R}$-Cartier divisor, 
and $D$ a Cartier divisor.  
If $D -(K_X + B)$ is $f$-ample, 
then $R^pf_*(I(X,B)\mathcal{O}_X(D)) = 0$ for $p > 0$.
\end{Thm}

The above theorem was first proved by Nadel \cite{Nadel} in a more 
general complex analytic setting (but on a smooth underlying variety).
Indeed the multiplier ideal sheaf is defined for a line bundle $L$ 
with a singular hermitian metric on a complex manifold $X$.
A {\em singular hermitian metric} $h$ is a degenerate hermitian 
metric which can be written 
locally as $h = h_0e^{-\phi}$, where $h_0$ is a $C^{\infty}$ metric 
and $\phi$ is a locally integrable weight function.
The {\em multiplier ideal sheaf} $I(L,h)$ is defined as
the largest ideal sheaf such that all local sections of 
$I(L,h)L$ satisfy locally the $L^2$ condition with respect to the metric.
It is proved to be a coherent sheaf on $X$.

For example, assume that $X$ is an open subset of $\mathbf{C}^n$ and 
$B = \sum b_iB_i$ is an effective $\mathbf{R}$-divisor.
If the prime divisors $B_i$ have local equations $g_i = 0$,
then we can define an {\em algebraically defined singular hermitian metric}
on a trivial bundle $L$
by using the weight function $\phi = \sum_i b_i \log \vert g_i \vert$.
In this case, the algebraic and analytic multiplier ideal sheaves 
coincide: $I(X,B) = I(L,h)$.

The point is that there are a lot more metrics which are essentially different
from algebraic metrics.
One can use a limit of a sequence of algebraic metrics to 
produce a non-algebraic one provided that we can prove certain convergence.
For example, the analytic Zariski decomposition (\cite{Tsuji-ZD}, \cite{DPS}) 
is a singular hermitian metric which is defined naturally for 
an arbitrary pseudo-effective line bundle.
The metric on a Hodge bundle is another example of analytic metrics.

\vskip 1pc

A new extension technique using the multiplier sheaves 
was developed by Siu \cite{Siu-Pm-gt}
when he proved the deformation invariance of plurigenera:

\begin{Thm}
Let $f: X \to T$ be a smooth projective morphism.
Then the plurigenus $\dim H^0(X_t, mK_{X_t})$ 
of a fiber $X_t = f^{-1}(t)$ is independent of $t \in T$ for any
positive integer $m$.
\end{Thm}

Nakayama \cite{Nakayama-Pm} proved that positive solutions for conjectures
in the MMP including the Abundance Conjecture imply the invariance
of plurigenera.
But Siu proved the theorem in one step without using an inductive 
approach of the MMP.
The theorem is proved in the case where a fiber $X_t$ is of general type 
in \cite{Siu-Pm-gt}, and in general in \cite{Siu-Pm-non-gt} (see also
\cite{Takayama-Pm} and \cite{Paun}).
The vanishing theorem used in these proofs is the 
Ohsawa-Takegoshi type extension theorem
(\cite{OT}, \cite{Siu-Pm-gt}):

\begin{Thm}[Ohsawa-Takegoshi Extension Theorem]\label{OTE}
Let $X$ be a Stein manifold, and $Y$ 
a smooth hypersurface defined by a 
bounded holomorphic function $t$.
Let $(L,h)$ be a line bundle on $X$ with a singular hermitian metric, 
and $s \in H^0(Y, K_Y+L \vert_Y)$ a holomorphic section.
Assume that the curvature current $- dd^c \log h$ is 
semipositive as a real current of type $(1,1)$, and that 
$\int_Y \vert s \vert^2 h < \infty$.
Then there exists an extension $\tilde s \in H^0(X,K_X+L)$ 
such that $\tilde s = s \wedge dt$ on $Y$ and
\[
\int_X \vert \tilde s \vert^2 h \le 
C \sup \vert t \vert^2 \int_Y \vert s \vert^2 h
\]
for a universal constant $C$.
\end{Thm}

The advantage of Theorem~\ref{OTE} is that we do not need to assume 
the strict positivity of the metric, 
which corresponds roughly to the ampleness or bigness.
But the semipositivity for the metric is strictly stronger the
corresponding algebraic concept of the nefness. 
This is the reason why there is still no algebraic proof for the 
invariance theorem of plurigenera.

\vskip 1pc

The idea of the proof of the invariance of plurigenera
in the case where the fibers are of general type is to use 
a construction downstairs which is similar to the Zariski decomposition.
The spaces of sections of pluricanonical systems on the central fiber
define a series of metrics on the canonical line bundle, 
hence a sequence of multiplier ideal sheaves.
Similar constructions on the total space of the deformations
define different metrics, hence a 
different sequence of multiplier ideal sheaves. 
A suitable vanishing theorem relates these ideals, 
and the extension of sections is proved.

This proof allowed an algebraic analogue in 
\cite{canonical} (see also \cite{Nakayama-ZD} and \cite{extension}) 
because the canonical line bundle is big in this case.
The algebraic version of the extension theorem were generalized 
to the logarithmic situation in \cite{HM-bounded} and 
\cite{Takayama-bounded} (see also \cite{Tsuji-bound}).
The logarithmic extension theorem was 
used in the proof of the PL flip theorem \cite{HM-flip} explained later.
We note that the algebraization of the proof in the general case, i.e., non
general type case, is still an open problem.

Now we state the logarithmic extension theorem:

\begin{Thm}\label{extension}
Let $f: X \to T$ be a projective morphism from a smooth variety to an affine 
variety, and 
$B = \sum b_iB_i$ a $\mathbf{Q}$-divisor whose support is a 
normal crossing divisor and such that only one coefficient $b_0$ is equal to 
$1$ and other coefficients satisfy $b_i \in (0,1)$ for $i \ne 0$.
Set $Y=B_0$.
Let $r$ be a positive integer such that $r(K_X+B)$ has integral coefficients.
Set $(K_X+B) \vert_Y = K_Y+B_Y$.
Assume the following conditions:

\begin{enumerate}

\item $K_X+B$ is pseudo-effective with respect to $(X,\ulcorner B \urcorner)$.

\item $B-Y$ is big with respect to $(X,Y)$.

\end{enumerate}

Then natural homomorphisms
\[
H^0(X,mr(K_X+B)) \to H^0(Y,mr(K_Y+B_Y))
\]
are surjective for all positive integers $m$.
\end{Thm}

This is a correct log generalization of the extension theorem,
while \cite{extension} Example~4.3 showed that a naive extension is false.


\section{Flip theorem}

We start to explain how the remaining two conjectures on the flips, 
the existence and the termination, 
are proved for the scaled MMP 
under the additional assumption that the boundary is big.

The basic idea is to use the induction on the dimensions using the 
{\em adjunction formula} which relates canonical divisors in different 
dimensions: if $Y$ is a smooth divisor on a smooth variety $X$, 
then we have
\[
(K_X + Y)\vert_Y = K_Y.
\]
This is a very different approach from the proof in dimension three 
in \cite{Mori-flip}.

We need two theorems, the Special Termination Theorem and the reduction to 
PL flips, due to Shokurov \cite{Shokurov-PL} 
(see also \cite{Fujino}) preceding to the proof of the flip theorem.
The first one is on the termination of flips along boundary components
with coefficient $1$:

\begin{Thm}[Special Termination]
Let $(X,B)$ be a $\mathbf{Q}$-factorial DLT pair of dimension $n$
which is projective over a base.
Assume that the MMP holds, i.e., the existence and termination conjectures 
hold, in dimension less than $n$.
Let 
\[
(X,B) = (X_0,B_0) \dashrightarrow (X_1,B_1) \dashrightarrow \cdots
\]
be an infinite sequence of flips.
Then there exists a positive integer $m_0$ such that the flip
$(X_m,B_m) \dashrightarrow (X_{m+1},B_{m+1})$ is an isomorphism in a
neighborhood of $\llcorner B_m \lrcorner$ for every $m \ge m_0$.

Similarly, if the MMP holds in dimension less than $n$ under the additional 
condition that the MMP is scaled or the boundary is big, then the special 
termination holds under the additional assumption that the MMP is scaled or
the boundary $B$ is big with respect to $(X,B)$.
\end{Thm}

The idea of the proof is to use the adjunction formula to reduced boundary
components.
It turned out that the Special Termination Theorem, 
not the general termination, is sufficient to prove the
existence of minimal models thanks to the Non-Vanishing Thoerem explained in
the next section.

\vskip 1pc

A small contraction morphism of a $\mathbf{Q}$-factorial DLT pair 
$\phi: (X, B) \to Y$ in the MMP
is said to be a {\em PL} (prelimiting) contraction, 
if there is an irreducible component $S$ of $\llcorner B \lrcorner$ 
such that $-S$ is $\phi$-ample.
A flip for a PL contraction is called a {\em PL flip}.
The second one is the reduction to PL flips:

\begin{Thm}
Assume that the special termination theorem for the scaled MMP with big
boundaries and 
the existence of the PL flip hold in dimension $n$.
Then the existence of the flip in general holds in the same dimension. 
\end{Thm}

The idea of the proof is to introduce additional artificial boundary 
as a scaling, 
and then reduce the boundary back to the original state gradually 
by using the scaled MMP.
This is the first appearance of the MMP with scaling.

Hacon and McKernan \cite{HM-flip} proved the existence of a PL flip:

\begin{Thm}
Assume that a minimal models exists for any KLT pair of dimension $n-1$ 
when the boundary is big. 
Then a flip exists for arbitrary PL small contraction in dimension $n$.
\end{Thm}

The idea of proof is as follows.
Let $\phi: (X, B) \to Y$ be a PL contraction morphism.
We may assume that $B$ is a $\mathbf{Q}$-divisor.
Let $S$ be an irreducible component of $\llcorner B \lrcorner$ 
such that $-S$ is $\phi$-ample.
We use the adjunction formula: $(K_X + B)\vert_S = K_S + B_S$,
where $B_S$ is usually larger than the restriction $B \vert_S$ 
because of the singularities of $X$, namely the {\em subadjunction}. 
In order to prove the finite generation of the relative canonical ring
\[
R(X/Y, K_X + B) = \bigoplus_{m=0}^{\infty} \phi_*\mathcal{O}_X
(\llcorner m(K_X + B) \lrcorner)
\]
it is sufficient to prove the finite generation of 
the image of the restriction homomorphism 
\[
R(X/Y, K_X + B) \to R(S/Y, K_S + B_S).
\]
The target of this ring homomorphism is finitely generated 
by the induction hypothesis, 
but the difficulty is that the homomorphism is {\em not} surjective.
The point of the proof of the PL flip theorem is to identify the image of this
homomorphism.

By restricting $m$-canonical systems on $X$ to $S$ and cancelling fixed 
components from the boundary $B_S$, we obtain a series of 
boundaries $B_{S',m}$ depending on $m$ on a fixed birational model $S'$ of $S$.
By using the extension theorem (Theorem~\ref{extension}), 
we prove that all the pluricanonical forms 
with respect the boundaries $B_{S',m}$ on $S'$ are 
in the images of the restriction homomorphisms.
By the Polytope Decomposition Theorem in dimension $n-1$, 
we find a fixed birational model $S''$ 
of $S$ which dominates all the minimal models corresponding to this series of 
boundary divisors.
Then we prove that certain stability of the movable parts of 
$m$-canonical systems holds
when $m$ goes to infinity by using the vanishing theorem
and the effective version of the Base Point Free Theorem by Koll\'ar 
\cite{Kollar-effective}.

We note that arbitrary divisor is big for a birational map 
such as $\phi$.
The difficulty concerning the infinity 
arising from the fact that we should consider all the 
$m$-canonical systems at the same time 
is solved by the finiteness statement of the 
Polytope Decomposition Theorem.


\section{Termination of flips and the non-vanishing theorem}

An old approach to the termination conjecture is to use invariants of 
singularities initiated by Shokurov \cite{Shokurov-NVT}.
Let us consider the simplest case; let $X$ be a three dimensional variety 
with only terminal singulsrities, and $\mu: Y \to X$ a resolution
of singularities.
We write $\mu^*K_X = K_Y + \sum e_jE_j$.
By the assumption, we have $e_j < 0$ for all $j$.
The {\em difficulty} of the variety $X$ is the number of $j$ such that
$-1 < e_j$.
One can prove that the difficulty is a well defined non-negative integer, and
it decreases strictly after a flip.
The termination follows immediately.

If the singularity is worse or if the dimension is higher, then 
similar properties, the well-definedness and the monotoneness, fail.
But we can modify the definition of the difficulty and obtain   
some termination theorems.
The following is the known results up to now obtained by using the 
concept of difficulty:

\begin{Thm}
(1) (\cite{termination}, \cite{Shokurov-3fold-model}) 
The termination conjecture holds for arbitrary DLT pairs in dimension three. 

(2) (\cite{AHK})
The termination conjecture holds for a KLT pair $(X,B)$ of dimension four if 
$- (K_X+B)$ is numerically equivalent to an effective $\mathbf{R}$-divisor
over the base.

(3) (\cite{AHK}) If $(X, \sum b_iB_i)$ is a KLT pair of dimension four
such that $c_0K_X+\sum c_iB_i$ is relatively big over the base for 
some numbers $c_i \in \mathbf{R}$.
Then there exists a process of the MMP with scaling which terminates.
\end{Thm}

In particular, 
the existence of a minimal model for arbitrary DLT pair $(X,B)$ over
a base is proved only in the case $\dim X = 3$ at the moment.

\vskip 1pc

Birkar, Cascini, Hacon and McKernan \cite{BCHM} took a very different
approach to the termination conjecture.
They did not prove the termination as an isolated statement, but rather
included it into a chain of statements concerning the MMP.
By using the induction on dimension, the termination is reduced to the 
special termination as in the following way.

Let $(X,B)$ be a DLT pair which is projective over a base $T$.
Assume that the log canonical divisor $K_X+B$ is numerically equivalent
to an effective $\mathbf{R}$-divisor $M$ over $T$.
Then the minimality question for $(X,B)$ is equivalent to that for $(X,B+tM)$ 
for $t > 0$ as long as the the new pair is LC.
Assume for simplicity that $(X,B+tM)$ is DLT for some $t > 0$ and that 
the support of $\llcorner B+tM \lrcorner$ coincides with that of $M$.
By applying the Special Termination Theorem, we conclude that the flips 
terminates near $\llcorner B+tM \lrcorner$.
Since $M$ is numerically equivalent to $K_X+B$, it follows that $K_X+B$ 
becomes nef in this process.
In the general case, we need more careful argument on the support of $M$
during the process of the MMP with scaling.

There is a modified version of the above termination argument 
in Birkar \cite{Birkar}.
In particular it is proved that, if $K_X+B$ is
numerically equivalent to an effective $\mathbf{R}$-divisor, then 
a minimal model exists if $\dim X \le 5$.

\vskip 1pc

The remaining problem is to prove the existence of an 
effective $\mathbf{R}$-divisor $M$.
This is a generalization of the Non-Vanishing Theorem:

\begin{Thm}
Let $(X,B)$ be a KLT pair which is projective over a base $T$.
Assume the following conditions:

\begin{enumerate}

\item $K_X+B$ is pseudo-effective over $T$.

\item $B$ is big over $T$.

\end{enumerate}

Then there exists an effective $\mathbf{R}$-divisor $M$
which is numerically equivalent to $K_X+B$ over $T$.
\end{Thm}

The assertion is only a part of the existence theorem of a minimal model,
but this is the key point in its proof.
The proof in \cite{BCHM} 
is also a generalization of the Non-Vanishing Theorem which is
a part of the Base Point Free Theorem which is already proved when $K_X+B$
is nef, together with the idea similar to
the proof of the flip theorem.
We use the minimal model theorem which is already proved 
under the additional assumption that $K_X+B$ is effective,
and the Polytope Decomposition Theorem in dimension one less.

The proof proceeds roughly as follows.
As in the case of the Non-Vanishing Thorem or the reduction theorem 
to the PL flip, we increase the boundary artificially so that the pair
$(X,B)$ becomes DLT and $\llcorner B \lrcorner = Z$ is irreducible.
We consider the case where $B$ is a $\mathbf{Q}$-divisor for simplicity.
If we add an ample $\mathbf{Q}$-divisor $\epsilon H$ 
to the boundary, then $K_X+B+\epsilon H$ becomes big, and there exists
a minimal model for the pair $(X,B+\epsilon H)$.
If we take the limit $\epsilon \to 0$, then there may be infinitely 
many chambers in the space of divisors corresponding to the minimal models
of the pairs $(X,B+\epsilon H)$.
By using the Special Termination Theorem, 
we can prove that certain neighborhood of $Z$
becomes stable under this infinite chain of wall crossings.
Then we can fix a minimal model of $Z$, on which there exists a desired 
section thanks to the usual Base Point Free Theorem.
By using the Vanishing Theorem, we infer that 
this section is extended to a birational model of $X$.

The assumption that the boundary $B$ is big is indispensable 
in the first step of the above proof where we find $Z$.
Indeed as is shown in \cite{Birkar}, the Non-Vanishing Thorem is the most 
difficult point if we try to extend the minimal model theory to the case
where the boundary is not necessarily big.


\section{Birational maps between birational models}

The output of the MMP, a minimal model or a Mori fiber space, 
is not uniquely determined in general when we start with a fixed pair,
because there are choices of the faces in the process of the MMP.
We would like to control the non-uniqueness.
The answer is given by the Polytope Decomposition Theorem again.

\vskip 1pc

The non-uniqueness of a minimal model for a fixed variety 
is a new phenomena in dimension three or higher.
The following theorem asserts that the existence of flops is the only reason.

A {\em flop} for a pair $(X,B)$ is a diagram
\[
\begin{CD}
X @>{\phi}>> Y @<{\phi^+}<< X^+
\end{CD}
\]
which is a flip for another pair $(X,B')$, 
where $B'$ is a suitably chosen different 
boundary, and such that $K_X+B$ is numerically trivial for $\phi$;
$[K_X+B] = 0$ in $N^1(X/Y)$. 

\begin{Thm}
Let $f: (X,B) \to T$ be a projective morphism from a  KLT pair, let
$g: (Y,C) \to T$ and $g': (Y',C') \to T$ be its minimal models, and let
$\alpha: Y \dashrightarrow Y'$ be the induced birational map over $T$.
Assume that there exists a canonical model of $(X,B)$.
Then $\alpha$ is decomposed into a sequence of flops in the follwoing way:
there exists an effective $\mathbf{Q}$-Cartier divisor $D$ on
$Y$ such that $(Y,C+D)$ is still KLT and such that $\alpha$ becomes a 
composition of a sequence of birational maps 
\[
Y = Y_0 \dashrightarrow Y_1 \dashrightarrow \cdots \dashrightarrow Y_l
= Y'
\]
such that $\alpha_k: Y_{k-1} \to Y_k$ ($1 \le k \le l$) is a flop for the 
pair $(Y_{k-1},C_{k-1})$ as well as a flip for the pair
$(Y_{k-1},C_{k-1}+D_{k-1})$, where $C_{k-1}$ and $D_{k-1}$ are 
the strict transforms of $C$ and $D$, respectively.
\end{Thm}

We note that $\alpha$ is an isomorphism in codimension one as we already know.
We remark that the boundary $B$ need not to be big. 
There is a differnt version of the factorization theorem in \cite{flop}, where
we do not need to assume the existence of a canonical model, but we have to
assume that $B$ is a $\mathbf{Q}$-divisor.

A {\em marked minimal model} is a pair consisting of a minimal model 
and a birational map to a fixed reference model.
The number of birationally equivalent marked minimal models is finite if 
the boundary $B$ is big, but it is not the case in general (\cite{CY}).
The above theorem claims that there are still only finitely many marked 
minimal models which lie between two different minimal models.
It is conjectured that the number of birationally equivalent minimal models 
is finite {\em up to isomorphisms}, i.e., when we forget the markings. 

We use the Polytope Decomposition Theorem to prove the above theorem 
in the following way.
We take a general ample $\mathbf{R}$-Cartier divisor $H$ and $H'$ 
on $Y$ and $Y'$, respectively, and
we consider a triangle spanned by $C$, $C+H$ and $C+H'$ in the 
space of divisors on $Y$, 
where the strict transform of $H'$ on $Y$ is denoted by 
the same letter.
The canonical model coresponds to the chamber $\{0\}$.
If we choose $H$ and $H'$ small enough, then the closures of 
the chambers $V$ and $V'$ 
corresponding to the models $Y$ and $Y'$ contain $0$.
Moreover the union of the chambers between $V$ and $V'$ 
which contain $0$ in the closures contains the line segment joining 
the points $C+H$ and $C+H'$.
Then the wall crossing process provides a decomposition of $\alpha$.
By induction on the relative Picard numbers, we obtain eventually the 
decomposition to flops.

\vskip 1pc

As for the Mori fiber spaces, the non-uniqueness phenomenon 
appears already in dimension two.
For example, ruled surfaces $\mathbf{P}(\mathcal{O}_{\mathbf{P}^1}
\oplus \mathcal{O}_{\mathbf{P}^1}(d))$ over $\mathbf{P}^1$ 
for different integers $d$ are all birationally equivalent.
An {\em elementary transformation} of a ruled surface is a combination of
a blowing up at a point in a fiber and the blowing down of the strcit 
transform of the fiber.
If $g_i: X_i \to C$ for $i=1,2$ 
are ruled surfaces over the same curve $C$, then 
they are connected each other by a sequence of elementary transformations.
There is a different kind of decompositions; 
any birational map 
$\mathbf{P}^2 \dashrightarrow \mathbf{P}^2$ is 
decomposed into linear and quadratic transformations.
The latter is further decomposed into point blowings up, point blowings down, 
and elementary tansformations.

The {\em Sarkisov program} is a higher dimensional generalization. 
We can decompose any birational map between the total spaces of Mori fiber 
spaces into {\em elementary links},
generalizations of elementary transformations.
 
Hacon and McKernan \cite{HM-Sarkisov} proved the following:

\begin{Thm}
Let $f: (X,B) \to T$ be a projective morphism from a  KLT pair, let
$\phi: (Y,C) \to Z$ and $\phi': (Y',C') \to Z'$ 
be Mori fiber spaces obtained 
by the MMP from a KLT pair $(X,B)$ over $T$, and let
$\alpha: Y \dashrightarrow Y'$ be the induced birational map over $T$:
\[
\begin{CD}
Y @>{\alpha}>> Y' \\
@V{\phi}VV @VV{\phi'}V \\
Z @. Z'
\end{CD}
\]
Then there exists a sequence of the following type commutative diagrams 
called {\em elementary links} for $1 \le k \le l$ with some positive integer 
$l$:
\[
\begin{CD}
U_1^{(k)} @>{\beta^{(k)}}>> U_2^{(k)} \\
@V{f_1^{(k)}}VV @VV{f_2^{(k)}}V \\
V_1^{(k)} @. V_2^{(k)} \\
@V{g_1^{(k)}}VV @VV{g_2^{(k)}}V \\
W^{(k)} @>=>> W^{(k)}
\end{CD}
\]
where $\beta^{(k)}$ is a composition of a sequence of flips for a suitably 
chosen boundary on $U_1^{(k)}$ and 
such that $\alpha$ is decomposed as 
\[
\alpha = \alpha_l \circ \cdots \circ \alpha_1
\]
where each $\alpha_k$ is a birational map described 
in one of the following cases:

\begin{enumerate}

\item $g_1^{(k)}$ and $f_2^{(k)}$ are Mori fiber spaces, $f_1^{(k)}$ 
is a divisorial contraction, $g_2^{(k)}$ is a morphism with relative
Picard number one,
and $\alpha_k = \beta^{(k)} \circ (f_1^{(k)})^{-1}$.

\item $g_1^{(k)}$ and $g_2^{(k)}$ are Mori fiber spaces, 
$f_1^{(k)}$ and $f_2^{(k)}$ are divisorial contractions, 
and $\alpha_k = f_2^{(k)} \circ \beta^{(k)} \circ (f_1^{(k)})^{-1}$.

\item $f_1^{(k)}$ and $g_2^{(k)}$ are Mori fiber spaces, 
$f_2^{(k)}$ is a divisorial contraction, 
$g_1^{(k)}$ is a morphism with relative Picard number one,
and $\alpha_k = f_2^{(k)} \circ \beta^{(k)}$.

\item $f_1^{(k)}$ and $f_2^{(k)}$ are Mori fiber spaces, 
$g_1^{(k)}$ and $g_2^{(k)}$ are morphisms with relative Picard number one,
and $\alpha_k = \beta^{(k)}$.

\end{enumerate}
\end{Thm}

The theorem was proved in dimension three by Corti \cite{Corti}.
The crucial point of the proof is to prove that a sequence of elementary links
terminates. 
\cite{HM-Sarkisov} interpreted the sequence of elementary links as 
a wall crossing process, 
and proved the termination by using the Polytope Decomposition Theorem.
It is remarkable that the four types of elementary links which are 
apparently different have the same 
interpretation as a wall crossing process in terms of 
the polytope decomposition.

The argument of \cite{HM-Sarkisov} is as follows.
We take a general ample $\mathbf{R}$-Cartier divisor $H$ (resp. $H'$) 
on $Y$ (resp. $Y'$) such that $K_Y+C+H = \phi^*L$ 
(resp. $K_{Y'}+C'+H'=(\phi')^*L'$) for some 
ample $\mathbf{R}$-Cartier divisor $L$ (resp. $L'$) on $Z$ (resp. $Z'$),
and let $D$ (resp. $D'$) be the strict transforms of $H$ (resp. $H'$) on $X$.
Then the pair $(Y,C+H)$ (resp. $(Y',C'+H')$) is a minimal model of $(X,B+D)$
(resp. $(X,B+D')$), and $Z$ (resp. $Z'$) is the canonical model.

We consider a triangle $\bar V$ spanned by $B$, $B+D$ and $B+D'$ in the 
space of divisors on $X$, 
and let $V$ be the subset corresponding to the pseudo-effective log canonical 
divisors.
The points $(X,B+D)$ and $(X,B+D')$ are on the boundary of $V$.
We consider a path on the boundary which connects these points, 
and look at the chambers whose closures intersect this path.
Then the wall crossing process provides a decomposition of $\alpha$.


Department of Mathematical Sciences, University of Tokyo, 

Komaba, Meguro, Tokyo, 153-8914, Japan 

kawamata@ms.u-tokyo.ac.jp

\end{document}